\documentclass{article}
\usepackage{amsmath}
\usepackage{latexcad}
\usepackage{amsfonts}
\usepackage{amssymb}
\usepackage{graphicx}

\begin{document}

\title{On partita doppia}
\author{Piergiulio Katis\\School of Mathematics and Statistics\\University of Sydney\\Sydney, NSW 2006 Australia\\
\and N. Sabadini\\Dipartimento di Scienze dell'Informazione\\Universit\`a di Milano\\via Comelico, 39\\30135 Milano, Italy\\
\and R. F. C. Walters\\School of Mathematics and Statistics\\University of Sydney\\Sydney, NSW 2006 Australia}
\date{14th September 1998}
\maketitle

\section{Introduction.}

\subsection{Double-entry bookkeeping}

In 1494 Fra Luca Pacioli published in Venice one of the first printed
mathematical books [P1494]. One section, Computis e Scripturis, is the first
published description of \emph{partita doppia} or double-entry bookkeeping,
the foundation of accounting. Double-entry bookkeeping had been developed over
a period of 200 years by Italian merchants and bankers. The aim of accounting
is the measurement of a distributed concurrent system, and it is our
contention that it is one of the earliest and most successful mathematical
theories of concurrency. It is interesting that at that time negative numbers
were not accepted; in the following century mathematicians repaired this
deficiency though there seems to have been no attempt to explicate the
distributed algebra underlying partita doppia. For further information about
early accounting see [BJ84], [DeR56].

\noindent In this paper we give a precise mathematical account of partita
doppia in terms of an algebraic structure on \textbf{Span(RGraph)} - the
bicategory of spans of reflexive graphs. Some interesting new mathematical
considerations arise in the study of this example. The fact that accounting
concerns the measurement of distributed systems is explicit in the form of the
mathematics. We need to consider \textbf{Span(RGraph)} as a not-necessarily
self-dual compact closed bicategory in order to account for the direction of
flow of value. Another new aspect is that accounts give an example of a
compositional concept represented by a lax functor: we are aware of many other
compositional concepts similarly represented and believe that this is a
widespread phenomenon, one which we will discuss in further detail in a
complete version of this paper. There is also an interesting new mathematical
result regarding the compact-closedness of a certain lax comma construction.

\noindent Another aspect to note in this work is the essential use of the two
dimensional structure of \textbf{Span(RGraph)}. Locally ordered bicategories
(whose arrows are `relations extended in time') have been used by Abramsky to
model interaction [Abr93], but in this paper a local order is insufficient:
the crucial totalling of values in a system of accounts is a genuine 2-cell.

\subsection{Pacioli}

Pacioli is an interesting figure, born in Borgo San Sepolcro (now Sansepolcro)
circa 1452, he was a student of another son of that city, Piero della
Francesca. He taught in Florence, Rome, Perugia, Milan, Venice and Bologna. He
was a friend of Leonardo, and in fact they collaborated on the book, De Divina
Proportione (Venezia, 1509). His \emph{Summa} was highly successful and
assisted the dispersion of mathematics, contributing in this way to the
flowering of mathematics in Europe in the 16th century. The book is a
compilation of known results. Pacioli is known as the father of accounting.

\noindent Vasari in \emph{The lives of the artists} (1568), was suspicious of
Pacioli's originality: \leftline{From the chapter of Vasari [V1568]on Piero
della Francesca:}

\leftskip=0.25in\rightskip=0.25in {\small
%
%
%
``The Franciscan, Luca dal Borgo, who wrote about the regular geometrical
bodies, was his [Piero della Francesca's] pupil; and when Piero died at an
advanced age, having written many books, this Luca arrogated them to himself
and published as his own work what had fallen into his hands after his
master's death.''}

\leftskip=0in\rightskip=0in

\noindent One section of De Divina Proportione seems to be essentially due to
Piero della Francesca.

\subsection{Prerequisites}

Graphs have been used since time immemorial to model systems, with the
vertices representing \emph{states} and the edges representing \emph{atomic
actions} of the system. Recently (in [KSW97b], [KSW97c]) the authors have
proposed algebraic structures on \textbf{Span(RGraph)} as a way of modelling
non-deterministic concurrent systems. We recall here briefly the ideas in this
model but refer for more details to those two papers.

\noindent An initial version of this work on partita doppia was described by
the third author at New Trends in Semantics, Bologna, 4-5 July, 1997.

\noindent This work arose out of earlier study of bicategories of processes
[KSW97a], [K96] which itself arose from the study of distributive categories
and imperative programming [W92a], [W92b], [KW93].

\noindent Much of the notation in this paper originates in category theory;
for this the reader may consult a text such as [M70] or [W92a].

\subsection{Acknowledgment}

\noindent We are grateful in particular for conversations with Henry Weld. We
acknowledge also the support for this project by the Australian Research
Council, Italian MURST 40\% and the Italian CNR.

\section{Span(RGraph)}

The algebra of this paper is the compact-closed symmetric-monoidal bicategory
\textbf{Span(RGraph)}. Although these notions are well-known in the
categorical literature (compact-closed [KL80], monoidal [M70], bicategory
[B67], Span [B67], reflexive graphs [Law89]) we will describe them briefly
here at the same time introducing a pictorial representation for expressions
of arrows.

\subsection{The algebra of graphs}

\noindent The algebra \textbf{Span(RGraph)} consists of three sorts of things,
objects, arrows and 2-cells. We will list the operations of the algebra as we
describe its elements.

\subsubsection{Objects}

The objects of \textbf{Span(RGraph)} are (finite) \emph{reflexive} graphs:
that is, each object $X$ consists of a set $X_{1}$ of \emph{edges}, a set
$X_{0}$ of \emph{vertices}, and three functions, the \emph{domain} function
$\mathrm{d}_{0}:X_{1}\to X_{0}$,(assigning the beginning vertex to an edge),
the \emph{codomain} function $\mathrm{d}_{1}:X_{1}\to X_{0}$, (assigning the
end vertex to an edge), and the \emph{null} function $\epsilon:X_{0}\to X_{1}%
$, (assigning a null loop to each vertex). Further there are axioms to be
satisfied, namely that $\mathrm{d}_{0}\epsilon=1_{X_{0}}$ and $\mathrm{d}%
_{1}\epsilon=1_{X_{0}}$. There is an obvious notion of morphism of reflexive
graphs. Since in this paper all graphs are reflexive we shall for brevity drop
the adjective reflexive from now on. Graphs may be represented in the usual
way as geometric figures with vertices and edges.

\subsubsection{Operations on objects}

There are two operations on objects which produce objects: (i) the products of
graphs $X\times Y$ is the graph whose vertices (edges) are pairs of vertices
(edges), one from $X$ and one from $Y$; the reverse graph $X^{-1}$ of a graph
$X$, which is actually just the graph $X$. It may seem strange to consider an
operation whose effect is the identity, but it will enable us to make
distinctions we need to make later on.

\subsubsection{Arrows}

\noindent Given graphs $X$ and $Y$ an arrow from $X$ to $Y$ consists of a
graph $R$ and two graph morphisms $\partial_{0} : R\to X, \partial_{1} : R \to
Y$. Such an arrow is called a \emph{span} of graphs. It is often denoted as follows:

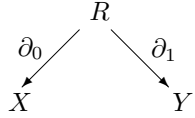
\centerline{
{\tt\setlength{\unitlength}{0.75pt}
\begin{picture}(110,79)
\thinlines\put(15,43){$\partial_0$}
\put(82,42){$\partial_1$}
\put(51,61){$R$}
\put(93,15){$Y$}
\put(10,15){$X$}
\put(62,56){\vector(1,-1){29}}
\put(46,56){\vector(-1,-1){29}}
\end{picture}}} \noindent We call the graph $R$ the head of the span and the
morphisms $\partial_{0}$ and $\partial_{1}$ the two legs of the span. We will,
by abuse of notation, denote the two legs of any span by the same symbols
$\partial_{0}$ and $\partial_{1}$. We will often denote the span simply as $R:
X\to Y$. We call $X$ and $Y$ the domain and codomain, respectively, of $R$; or
the boundaries of $R$.

\noindent We represent a span $R: X\to Y$ by a picture of the following form:

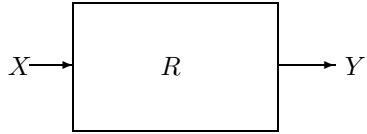
\centerline{
{\tt\setlength{\unitlength}{0.75pt}
\begin{picture}(197,84)
\thinlines\put(87,38){$R$}
\put(180,38){$Y$}
\put(10,38){$X$}
\put(147,43){\vector(1,0){28}}
\put(21,43){\vector(1,0){22}}
\put(43,10){\framebox(103,64){}}
\end{picture}}} \noindent We may also picture a span in another way. If the
objects are given as products of graphs, for example $X=X_{1}\times X_{2}%
^{-1}$, $Y=Y_{1}\times Y_{2}^{-1}\times Y_{3}$ we picture the span
%
%
as:

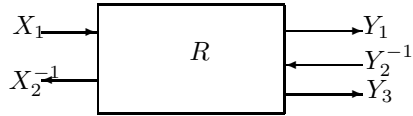
\centerline{
{\tt\setlength{\unitlength}{0.65pt}
\begin{picture}(240,83)
\thinlines\put(169,21){\vector(1,0){45}}
\put(214,38){\vector(-1,0){45}}
\put(169,58){\vector(1,0){45}}
\put(113,41){$R$}
\put(216,19){$Y_3$}
\put(215,36){$Y_2^{-1}$}
\put(214,56){$Y_1$}
\put(08,24){$X_2^{-1}$}
\put(10,56){$X_1$}
\put(59,29){\vector(-1,0){32}}
\put(27,57){\vector(1,0){32}}
\put(60,10){\framebox(108,63){}}
\end{picture}}} \noindent In this case we may call each of $X_{1}$,
$X_{2}^{-1}$, $Y_{1}$, $Y_{2}^{-1}$, $Y_{3}$, boundaries of $R$.

\subsubsection{\textbf{Operations involving objects and arrows}}

\subsubsection{\textbf{Composition of spans.}}

\noindent The composite of spans $R:X\to Y$ and $S:Y\to Z$ is the span
$R\bullet S : X\to Z$ whose head is the graph with vertex set
\[
\{(r,s); r\text{ is a vertex of } R, s\text{ is a vertex of }S, \partial
_{1}(r) = \partial_{0}(s)\}
\]
and with edge set
\[
\{(\rho,\sigma); \rho\text{ is a edge of }R, \sigma\text{ is a edge of } S,
\partial_{1}(\rho) = \partial_{0}(\sigma)\}.
\]
Beginnings and ends of edges have the obvious definitions. ($R\bullet S$ is
the pullback $R \times_{Y} S$).

\noindent The pictorial representation of the composition of two spans $R$ and
$S$ is as follows (with the obvious modification if the objects are products
of graphs!):

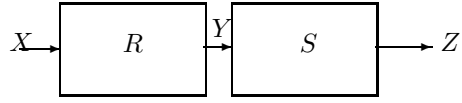
\centerline{{\tt\setlength{\unitlength}{0.65pt}
\begin{picture}(278,73)
\thinlines\put(40,10){\framebox(84,53){}}
\put(140,10){\framebox(84,53){}}
\put(16,37){\vector(1,0){24}}
\put(123,38){\vector(1,0){17}}
\put(223,38){\vector(1,0){33}}
\put(10,36){$X$}
\put(128,44){$Y$}
\put(261,36){$Z$}
\put(76,35){$R$}
\put(179,35){$S$}
\end{picture}}}

\smallskip\noindent\textbf{Tensor of spans.}

\noindent The tensor of two spans $R:X\to Y$ and $S:Z\to W$ is the span
denoted $R \otimes S : X \times Z\to Y \times W$, and defined by: the head of
$R \otimes S$ is $R \times S$; the legs of $R \otimes S$ are $\partial_{0}
\times\partial_{0}$ and $\partial_{1} \times\partial_{1}$. The pictorial
representation of the tensor of two spans is:

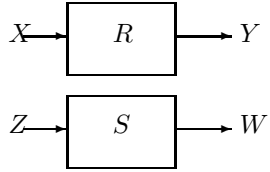
\centerline{{\tt\setlength{\unitlength}{0.50pt}
\begin{picture}(207,144)
\thinlines\put(94,36){$S$}
\put(94,105){$R$}
\put(190,36){$W$}
\put(15,36){$Z$}
\put(190,105){$Y$}
\put(15,105){$X$}
\put(141,39){\vector(1,0){43}}
\put(141,109){\vector(1,0){43}}
\put(26,39){\vector(1,0){33}}
\put(26,109){\vector(1,0){33}}
\put(60,10){\framebox(81,54){}}
\put(60,80){\framebox(81,54){}}
\end{picture}}
} \noindent In addition to these operations there are the following constants
of the algebra.

\smallskip\noindent\textbf{The spans $\eta_{X} : I \to X^{-1} \times X$ and
$\varepsilon_{X} : X \times X^{-1} \to I$.}

\noindent The terminal graph, denoted $I$, has one vertex, $0$ and one edge,
by necessity the null loop. The span with head $X$ and legs $!:X\to I,
\Delta:X\to X\times X$ is called $\eta_{X}$.

\noindent The span with head $X$ and legs $\Delta:X\to X\times X, !:X\to I $
is called $\varepsilon_{X}$.

\noindent The two spans are pictured thus:

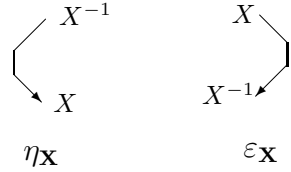
\centerline{
{\tt\setlength{\unitlength}{0.75pt}
\begin{picture}(187,100)
\put(122,11){\begin{picture}(55,79)
\thinlines\put(-15,30){$X^{-1}$}
\put(0,71){$X$}
\put(6,3){{\large$\mathbf{\varepsilon_X}$}}
\put(28,61){\line(0,-1){12}}
\put(28,50){\vector(-1,-1){16}}
\put(14,75){\line(1,-1){14}}
\end{picture}}
\thinlines\put(32,37){$X$}
\put(35,81){$X^{-1}$}
\put(17,13){\large{$\mathbf{\eta_X}$}}
\put(12,56){\vector(1,-1){14}}
\put(12,68){\line(0,-1){12}}
\put(28,84){\line(-1,-1){16}}
\end{picture}}}

\noindent The arrows $\eta$ and $\varepsilon$ are the unit and counit of the
compact-closed structure on \textbf{Span(RGraph)}. It will become clear that
their role in the context of this paper is to permit a feedback operation on
distributed systems.

\noindent The correspondence between constants and operations, and the
geometric representations given above, result in the fact that expressions in
the algebra have corresponding circuit or system diagrams. We will draw such
pictures later when discussing systems of accounts, though in this paper the
pictures are solely an aid to reason and not precisely formalized. The formal
thing is the expression.

\subsubsection{2-cells of Span(RGraph)}

\noindent If $R$ and $S$ are spans from $X$ to $Y$, a 2-cell $\phi:R\to S$
consists of a graph morphism $\phi$ between the heads of the spans satisfying
$\partial_{0} \phi= \partial_{0}$, and $\partial_{1} \phi=\partial_{1}$.

\noindent Corresponding to each span $R$ there is an obvious identity 2-cell
$R\rightarrow R$. Further 2-cells compose in two ways, horizontally and
vertically. Vertically, the composite $\phi\cdot\psi$ of $\phi:R\rightarrow S$
with $\psi:S\rightarrow T$ is formed by the composition of the graph morphisms
$\phi$ and $\psi$. Horizontally, the composite $\phi\bullet\psi$ of
$\phi:R\rightarrow S:X\rightarrow Y$ with $\psi:T\rightarrow U:Y\rightarrow Z$
is formed by pullback.

\subsection{Spans from graph morphisms}

\noindent Given a graph morphism $f:X\rightarrow Y$ there are two special
associated spans $f_{\ast}:X\rightarrow Y$ and $f^{\ast}:Y\rightarrow X$
defined as follows%

\begin{align*}
f_{\ast}  &  =(1_{X},f):X\rightarrow Y,\\
f^{\ast}  &  =(f,1_{Y}):Y\rightarrow X.
\end{align*}

\noindent and two 2-cells $1_{X}\rightarrow f_{\ast}\bullet f^{\ast}$,
$f^{\ast}\bullet f_{\ast}\rightarrow1_{Y}$ exhibiting $f^{\ast}$ as the right
adjoint of $f_{\ast}$. We will need the following straightforward proposition
in section 4.

\smallskip\noindent\textbf{Proposition}

\noindent If $f:U\rightarrow X$ and $g:V\rightarrow Y$ are graph morphisms and
$R:U\rightarrow V$, $S:X\rightarrow Y$ are spans of graphs then there is a
bijection between 2-cell of spans
\[
\phi:R\bullet g_{\ast}\rightarrow f_{\ast}\bullet S
\]

\noindent and graph morphisms
\[
\varphi:R\rightarrow S\text{ such that }\partial_{0}\varphi=f_{\ast}%
\partial_{0}\text{ and }\partial_{1}\varphi=g_{\ast}\partial_{1}.
\]

\subsection{Behaviours of a span}

\smallskip\noindent\textbf{Definition}

\noindent A behaviour $\pi$ of a span $R : X \to Y$ is a finite path in the
graph $R$, the head of the span.

\noindent Notice that applying the legs of the span to a behaviour $\pi$
yields two paths, one $\partial_{0}(\pi)$ in $X$ and the other $\partial
_{1}(\pi) $ in $Y$. The graphs $X$ and $Y$ may be thought of as the (left and
right) boundaries of the system $R$. Then $\partial_{0}(\pi)$, $\partial
_{1}(\pi)$ may be thought of as the behaviour of the boundaries of the system
corresponding to the behaviour $\pi$, or the behaviour of the system reflected
on the boundaries.

\noindent The following result is straightforward.

\smallskip\noindent\textbf{Proposition}

(i) A behaviour of the composite $R;S$ of two spans is a pair of behaviours,
one $\rho$ of $R$, the other $\sigma$ of $S$, such that $\partial_{1}(\rho) =
\partial_{0}(\sigma)$. That is, a behaviour of $R;S$ consists of a behaviour
of $R$ and a behaviour of $S$ which agree (synchronize) on the common boundary.

(ii) A behaviour of the tensor $R \otimes S$ of two spans is just a pair of
behaviours, one $\rho$ of $R$, the other $\sigma$ of $S$.

(iii) A behaviour of $\eta_{X} : I \to X \times X$ is a path in $X$ reflected
(synchronously and equally) on the two boundaries. The behaviours of
$\varepsilon$ are similarly described.

\section{Standard Accounts}

\noindent We need first to define the notion of standard account (in which a
general account will be valued).

\smallskip\noindent\textbf{Definition} A channel is a graph with one vertex
and edges being all non-negative integers. Suppose $X_{1}$, $X_{2}$, $\dots$,
$X_{m}$, $Y_{1}$, $Y_{2}$, $\dots$, $Y_{n}$ are channels, and suppose that
$X=X_{1}^{\xi_{1}}\times\dots X_{m}^{\xi_{m}}$ and $Y=Y_{1}^{\zeta_{1}}%
\times\dots Y_{n}^{\zeta_{n}}$ where the $\xi= (\xi_{1},\ldots,\xi_{m})$ and
$\zeta= (\zeta_{1}, \ldots, \zeta_{n})$ are sequences of $1$'s and $-1$'s.
Then the \emph{standard account} $A_{X,Y}$ from $X$ to $Y$ (which is often
denoted merely by $A$) is the span whose head vertex set is the integers, and
an edge $\rho: r\to s$ is an $m+n$-tuple of natural numbers $x_{i}
(i=1,2,\ldots,m)$, $y_{j} (j=1,2,\ldots,n)$ each satisfying
\[
s-r= \sum_{i=1}^{m} \xi_{i} x_{i} - \sum_{j=1}^{n} \zeta_{j} y_{j},
\]
and $\partial_{0,i}(\rho)=x_{i}$, $\partial_{1,j}(\rho)=y_{j}$ where
$\partial_{k,l}$ is $\partial_{k}$ followed by the $l$'th projection. \smallskip

\noindent In words, a vertex of an account is a possible \emph{value}, and an
edge is a change of value of the account as a result of various ingoings and
outgoings of value. The condition says that the change in the value of the
account after a transaction is the result of the difference between ingoings
and outgoings - that is, the condition is a \emph{continuity equation} for
value. Notice that value can only flow \emph{into} the accounts on the
channels of the form $X^{+1}$ on the left, and on channels of the form
$X^{-1}$ on the right, and \emph{out of} the accounts on channels of the form
$X^{-1}$ on the left, and on channels of the form $X^{+1}$ on the right.

\smallskip

\noindent\textbf{Note} Standard accounts do not form a sub-algebra of
\textbf{Span(RGraph)}; in particular, the collection of standard accounts does
not contain identities and it is not closed under composition nor tensor.
However, the class of standard accounts bears extra structure which will be
used in the next section to define a compact closed bicategory of accounts.
This structure is presented below in the form of data and axioms.

\smallskip

Before starting, note that if $\beta, \beta^{\prime}: G \rightarrow A_{X,Y}$
are 2-cells from any span of graphs to an account such that $\beta$ and
$\beta^{\prime}$ take the same value on the vertices of $G$, then $\beta=
\beta^{\prime}$. In the following, we will define families of 2-cells in
\textbf{Span(RGraph)} by only specifying their value on vertices. The reader
can easily verify the existence of these 2-cells.

\smallskip

\noindent\textbf{Data}

\begin{enumerate}
\item  For each account $A:X\rightarrow X$, let $\theta:1_{X}\rightarrow A$ be
the (unique) 2-cell which maps the only vertex of the head of $1_{X}$ to the
vertex $0\in A$.

\item  For any pair of accounts $A_{X,Y}$ and $A_{Y,Z}$, let $\alpha
:A_{X,Y}\bullet A_{Y,Z}\rightarrow A_{X,Z}$ be the 2-cell such that
$\alpha(i,j)=i+j$.

\item  For any pair of accounts $A_{W,X}$ and $A_{Y,Z}$, let $\tau
:A_{W,X}\otimes A_{Y,Z}\rightarrow A_{W\otimes Y,X\otimes Z}$ be the 2-cell
such that $\tau(i,j)=i+j$.

\item  For each account $A:I\rightarrow X^{-1}\otimes X$, let $\delta:\eta
_{X}\rightarrow A$ be the 2-cell which maps the only vertex of the head of
$\eta_{X}$ to the vertex $0\in A$. (Of course, if $X=X_{1}^{\xi_{1}}%
\times\dots X_{m}^{\xi_{m}}$ then $X^{-1}=X_{m}^{-1\times\xi_{m}}\times\dots
X_{1}^{-1\times\xi_{1}}$.)

\item  For each account $A:X\otimes X^{-1}\rightarrow I$, let $\gamma
:\epsilon_{X}\rightarrow A$ be the 2-cell which maps the only vertex of the
head of $\epsilon_{X}$ to the vertex $0\in A$.
\end{enumerate}

\smallskip

\noindent\textbf{Axioms}

\begin{enumerate}
\item  For any $X$ and $Y$, the following 2-cells are equal.
\[
(\theta\bullet A_{X,Y})\cdot\alpha:(1_{X}\bullet A_{X,Y})\rightarrow
A_{X,X}\bullet A_{X,Y}\rightarrow A_{X,Y}%
\]%
\[
1_{A_{X,Y}}:A_{X,Y}\rightarrow A_{X,Y}%
\]

\item  For any $W$, $X$, $Y$ and $Z$, the following 2-cells are equal.
\[
(\alpha\bullet A_{Y,Z})\bullet\alpha:(A_{W,X}\bullet A_{X,Y})\bullet
A_{Y,Z}\rightarrow A_{W,Y}\bullet A_{Y,Z}\rightarrow A_{W,Z}%
\]%
\[
(A_{W,X}\bullet\alpha)\cdot\alpha:A_{W,X}\bullet(A_{X,Y}\bullet A_{Y,Z}%
)\rightarrow A_{W,X}\bullet A_{X,Z}\rightarrow A_{W,Z}%
\]

\item  Given standard accounts $X\rightarrow Y$, $Y\rightarrow Z$, $X^{\prime
}\rightarrow Y^{\prime}$ and $Y^{\prime}\rightarrow Z^{\prime}$, the following
2-cells are equal.
\[
(\alpha\otimes\alpha)\cdot\tau:(A\bullet A)\otimes(A\bullet A)\rightarrow
A\otimes A\rightarrow A
\]%
\[
(\tau\bullet\tau)\cdot\alpha:(A\bullet A)\otimes(A\bullet A)\cong(A\otimes
A)\bullet(A\otimes A)\rightarrow A\bullet A\rightarrow A
\]

\item  For any $X$, the following 2-cells are equal:
\begin{align*}
{ ((\theta\otimes\delta)\bullet(\gamma\otimes\theta))\cdot(\tau\bullet
\tau)\cdot\alpha}  &  :{ (1}_{X}{ \otimes\eta}_{X}{ )\bullet(\epsilon}_{X}{
\otimes1}_{X}{ )}\\
&  \rightarrow{ (A\otimes A)\bullet(A\otimes A)\rightarrow A\bullet
A\rightarrow A}%
\end{align*}%
\[
\theta:(1_{X}\otimes\eta_{X})\bullet(\epsilon_{X}\otimes1_{X})=1_{X}%
\rightarrow A
\]

\item  For any $X$, the following 2-cells are equal:
\begin{align*}
{ ((\delta\otimes\theta)\bullet(\theta\otimes\gamma))\cdot(\tau\bullet
\tau)\cdot\alpha}  &  :{ (\eta}_{X}{ \otimes1}_{X}{ )\bullet(1}_{X}{
\otimes\epsilon}_{X}{ )}\\
&  \rightarrow{ (A\otimes A)\bullet(A\otimes A)\rightarrow A\bullet
A\rightarrow A}%
\end{align*}%
\[
{ \theta:(\eta}_{X}{ \otimes1}_{X}{ )\bullet(1}_{X}{ \otimes\epsilon}_{X}{
)=1}_{X^{-1}}{ \rightarrow A}%
\]
\end{enumerate}

\noindent\textbf{Remark} The structure described above is essentially that of
a lax morphism of compact closed bicategories. (The domain of this morphism is
the chaotic category whose objects are strings of plus and minus signs - the
compact closed structure being the obvious one. The codomain is, of course,
\textbf{Span(RGraph)}.) The first two data (and axioms) are part of the
structure of a morphism of bicategories. The second half of the data (and
axioms) relate to the compact closed structure. This abstract structure will
be defined and investigated in detail in another paper.

\section{The compact closed bicategory of accounts}

\noindent In this section we a make a formal definition of the compact-closed
bicategory \textbf{Accounts} of \emph{accounts}, in such a way that a
\emph{system of accounts with partita doppia} will be an expression in
\textbf{Accounts}. The connection with conventional accounting will be made in
the next section.

\smallskip\noindent\textbf{Definition}

\begin{itemize}
\item  An object of \textbf{Accounts} is a graph $U$ together with a
(reflexive) graph morphism $f:U\rightarrow X$ where $X$ is a products of
channels $X=X_{1}^{\xi_{1}}\times\dots\times X_{m}^{\xi_{m}}$.

\item  An arrow of \textbf{Accounts}, called a \emph{general account}, from
$f:U\rightarrow X$ to $g:V\rightarrow Y$ consists of a span $R:U\rightarrow
V$, and a 2-cell $\phi_{R}:R\bullet g_{\ast}\rightarrow f_{\ast}\bullet
A_{X,Y}$.
\end{itemize}

\centerline{\setlength{\unitlength}{0.85mm}
\begin{picture}(44,44)
\thinlines
\drawcenteredtext{9.0}{34.0}{$U$}
\drawcenteredtext{35.0}{34.0}{$V$}
\drawcenteredtext{9.0}{10.0}{$X$}
\drawcenteredtext{35.0}{10.0}{$Y$}
\drawcenteredtext{21.0}{6.0}{$A_{X,Y}$}
\drawcenteredtext{21.0}{38.0}{$R$}
\drawcenteredtext{5.0}{22.0}{$f_*$}
\drawcenteredtext{39.0}{22.0}{$g_*$}
\drawcenteredtext{19.15}{24.69}{$\phi_R$}
\drawvector{9.0}{32.0}{20.0}{0}{-1}
\drawvector{35.0}{32.0}{20.0}{0}{-1}
\drawvector{11.0}{34.0}{22.0}{1}{0}
\drawvector{11.0}{10.0}{22.0}{1}{0}
\drawpath{23.62}{23.44}{19.62}{19.44}
\drawpath{23.0}{24.0}{19.0}{20.0}
\drawpath{19.37}{21.66}{18.69}{18.97}
\drawpath{18.69}{18.97}{20.94}{19.65}
\end{picture}}

\begin{itemize}
\item  The composite of general accounts $(R:U\rightarrow V,\phi_{R})$,
$(S:V\rightarrow W,\phi_{S})$ is
\[
\big(R\bullet S:U\rightarrow W,(R\bullet\phi_{S})\cdot(\phi_{R}\bullet
A)\cdot(f_{\ast}\bullet\alpha):R\bullet S\bullet h_{\ast}\rightarrow f_{\ast
}\bullet A\big).
\]
\end{itemize}

\smallskip\noindent\textbf{Theorem} The bicategory \textbf{Accounts} is
compact closed.

\smallskip\noindent\textbf{proof} The proof (details will be given elsewhere)
amounts to first checking that left adjoint arrows in a compact closed
bicategory \textbf{B} are the objects of a compact closed bicategory
Ladj(\textbf{B}), the arrows being squares in the bicategory containing a
2-cell. Combining this with the lax structure on standard accounts yields the
required structure on \textbf{Accounts}.

\smallskip\noindent\textbf{Remark} As a result of this theorem we may draw
pictures of expressions of general accounts similar to those we described
earlier for expressions in \linebreak \textbf{Span(RGraph)}.

\smallskip\noindent\textbf{Definition} A \emph{system of accounts} is an
expression in the compact closed bicategory \textbf{Accounts}. A \emph{closed}
system of accounts is an expression with domain and codomain of the form
$!_{U}:U\to I$.

\smallskip\noindent\textbf{Corollary} For a closed system of accounts there is
an invariant of a behaviour, called the \emph{total value}.

\smallskip\noindent\textbf{proof} A closed system of accounts evaluates as a
span $R:U\rightarrow V$ and a 2-cell $\phi_{R}:R\bullet(!_{U})_{\ast
}\rightarrow(!_{V})_{\ast}\bullet A_{I,I}$. The 2-cell amounts to a graph
morphism from $R$ to $A_{I,I}$. But $A_{I,I}$ is a discrete reflexive graph,
and so each path in $R$ lies over a single vertex of $A_{I,I}$.

The fact that at this point the existence of an invariant is essentially
trivial is a consequence (and an evidence for the naturality of) the structure
of standard accounts.

\section{The relation with conventional bookkeeping}

\noindent A general account $(R,\phi_{R})$ is that it is a transition system
which is \emph{measured} by a standard account, the measure being the 2-cell
$\phi_{R}$. The state of a conventional account may be much more than just its
value; for example, an account may contain a record of its history. Further a
transaction involving two accounts usually has much more information than just
the value passed; for example, it may include the addresses of the people
involved in the transaction.

A conventional accounting system is a closed system consisting of an
expression in five different types of accounts - asset accounts, liability
accounts, equity accounts, expense and income accounts, with the following schematic

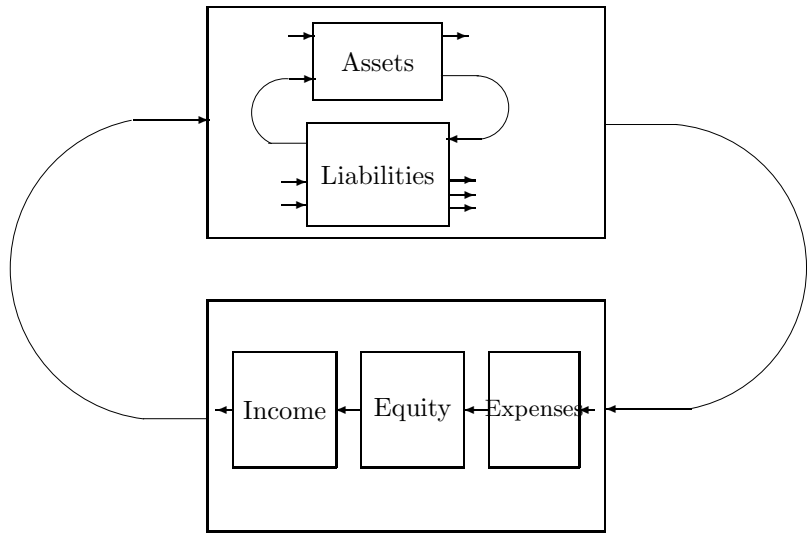
\centerline{\setlength{\unitlength}{0.85mm}
\begin{picture}(105,90)
\thinlines
\drawframebox{52.0}{68.0}{62.0}{36.0}{}
\drawframebox{52.0}{22.0}{62.0}{36.0}{}
\drawvector{9.3}{68.3}{11.9}{1}{0}
\drawvector{96.59}{23.13}{13.51}{-1}{0}
\drawpath{21.0}{21.56}{10.87}{21.56}
\drawpath{83.09}{67.63}{94.34}{67.63}
\drawarc{13.79}{45.08}{47.38}{1.69}{4.51}
\drawarc{92.09}{45.23}{44.99}{4.8}{1.37}
\drawframebox{47.54}{77.51}{20.0}{12.0}{Assets}
\drawframebox{47.55}{59.77}{22.0}{16.0}{Liabilities}
\drawvector{33.61}{81.56}{4.0}{1}{0}
\drawvector{33.73}{74.83}{4.0}{1}{0}
\drawvector{57.69}{81.56}{4.0}{1}{0}
\drawpath{57.45}{75.27}{63.3}{75.27}
\drawvector{32.54}{58.65}{4.0}{1}{0}
\drawvector{32.5}{55.04}{4.0}{1}{0}
\drawvector{58.81}{58.86}{4.0}{1}{0}
\drawvector{58.81}{56.61}{4.0}{1}{0}
\drawvector{58.81}{54.59}{4.0}{1}{0}
\drawframebox{33.0}{23.0}{16.0}{18.0}{Income}
\drawframebox{53.0}{23.0}{16.0}{18.0}{Equity}
\drawframebox{72.0}{23.0}{14.0}{18.0}{\small{Expenses}}
\drawvector{44.95}{22.92}{4.03}{-1}{0}
\drawvector{64.97}{22.92}{4.05}{-1}{0}
\drawvector{24.92}{22.92}{2.71}{-1}{0}
\drawvector{81.4}{22.92}{2.48}{-1}{0}
\drawpath{36.52}{64.7}{30.45}{64.7}
\drawvector{63.75}{65.38}{5.41}{-1}{0}
\drawarc{33.15}{69.52}{10.63}{2.0}{4.79}
\drawarc{63.06}{70.31}{9.93}{4.71}{1.43}
\end{picture}}

\noindent An \emph{asset} account is one in which the states have values in
the non-negative integers, called debits; a \emph{liability} account is one in
which the states have values in the non-positive integers, called credits.
Liability accounts generally record borrowings from external sources which
need to be repaid. \emph{Expense,} \emph{income} and \emph{equity} accounts
exist in order to produce a closed system. Expense and income accounts record
outgoings and\ ingoings from external sources which are free of any further
obligation. Expense accounts transactions reduce assets and increase expenses,
and hence expense accounts have debit state; income account transactions
increase assets and reduce income and hence income accounts have credit
states. In a transaction involving the ingoing of income an asset account
becomes more positive while the income account becomes correspondingly more
negative. The composite of the three (types of) accounts, income, expense and
equity may be called owner's equity and represents (with negative sign) the
total value of the business. At least once a year the income and expense
accounts are zeroized - that is, value is transferred from the expense
accounts to the equity accounts, and from the equity accounts to the income
accounts, thereby placing all the owner's equity in the equity accounts.

The behaviours considered have invariant (total) 0. When the expense and
income accounts have been zeroized (or composed with the equity account) this
means the invariant is expressed by the conventional equation

\centerline{\it Assets=Liabilities+Owner's Equity}

\noindent In this equation Liabilities and Owner's Equity the negative of the
totals in the liability and equity accounts, respectively. The Balance Sheet
contains the values of the states of the asset, liability, and owner's equity accounts.

\noindent Notice that it is crucial in defining a system of accounts that the
diagonal operation (a useful operation in describing other concurrent systems)
must be avoided. Copying money is counterfeiting and leads to a failure of the invariant.

\bigskip

\noindent\textbf{An Example}

\noindent Consider a system of accounts with one account of each type
discussed. Suppose initally the assets are 1000, and hence the equity is
-1000. Then consider the following sequence of tranactions - first a purchase
of 2000 value of consumable goods is made incurring a liability of -2000. Next
there is income of 1500, as a result of which the asset account increases by
1500, and the income account reduces by the same amount. Next 1000 of assets
are used to reduce the liability. Next the expenses are zeroized, and then the
same for the income. At the end of this sequence of transactions there are
1500 debit of assets, 1000 credit of liabilities, and an owner's equity of 500 credit.

\bigskip%

\begin{gather*}
\text{Asset, Liability, Expense, Income, Equity}\\
(1000,0,0,0,-1000)\\
\downarrow\\
(1000,-2000,2000,0,-1000)\\
\downarrow\\
(2500,-2000,2000,-1500,-1000)\\
\downarrow\\
(1500,-1000,2000,-1500,-1000)\\
\downarrow\\
(1500,-1000,0,-1500,1000)\\
\downarrow\\
(1500,-1000,0,0,-500)
\end{gather*}

\section{Bibliography}

{\small {\ \noindent[Abr93] S. Abramsky, Interaction categories (extended
abstract), in Theory and Formal Methods Workshop, Springer Verlag, 1993. }}

{\small \noindent[A94] A. Arnold, Finite transition systems, Prentice Hall,
1994. }

{\small \noindent[B67] J. B\`enabou, Introduction to bicategories, Reports of
the Midwest Category Seminar, Lecture Notes in Mathematics 47, pages 1--77,
Springer-Verlag, 1967. }

{\small \noindent[BJ84] R. Gene Brown, Kenneth S. Johnston, \emph{Paciolo on
accounting}, Garland Publishing, NY and London, 1984. }

{\small \noindent[CW87] A. Carboni and R.F.C. Walters, Cartesian Bicategories
I, Journal of Pure and Applied Algebra, 49, pages 11-32, 1987. }

{\small \noindent[DeR56] R. De Roover, The development of accounting prior to
Luca Pacioli according to the account-books of medieval merchants, pp114-184,
A.C. Littleton and B.S. Yamey, \emph{Studies in the history of accounting},
London, Sweet and Maxwell, 1956. }

{\small \noindent[KSW97a] P. Katis, N. Sabadini, R.F.C. Walters, Bicategories
of processes, Journal of Pure and Applied Algebra, 115, no.2, pp 141 - 178,
1997 }

{\small \noindent[KSW97b] P. Katis, N. Sabadini, R.F.C. Walters, Span(Graph):
A categorical algebra of transition systems, Proceedings, Algebraic
Methodology and Software Technology, SLNCS 1349, 307-321, 1997. }

{\small \noindent[KSW97c] P. Katis, N. Sabadini, R.F.C. Walters, Representing
P/T nets in Span(Graph), Proceedings, Algebraic Methodology and Software
Technology, SLNCS 1349, 322-336, 1997. }

{\small \noindent[K96] P. Katis, Categories and bicategories of processes, PhD
Thesis, University of Sydney, 1996. }

{\small \noindent[KL80] G.M. Kelly, M.L. Laplaza, Coherence for compact closed
categories, Journal of Pure and Applied Algebra, 19:193-213, 1980. }

{\small \noindent[KW93] W. Khalil and R.F.C. Walters, An imperative language
based on distributive categories II, Informatique Th\`eorique et Applications,
27, 503-522, 1993. }

{\small \noindent[Law89] F. W. Lawvere, Qualitative Distinctions between some
Toposes of Generalized Graphs, Proceedings of AMS Boulder 1987, Symposium on
Categories in Computer Science and Logic, Contemporary Mathematics 92 (1989),
261-299. }

{\small \noindent[M70] S. Mac Lane, Categories for the working mathematician,
Springer Verlag, 1970. }

{\small \noindent\lbrack No84] Ch. Nobes, Ed., \emph{The development of double
entry}, Garland Publishing, NY and London, 1984. }

{\small \noindent[P1494] Fra Luca Pacioli, \emph{Summa de arithmetica
geometrie proportioni e proportionalita}, Venezia, 1494. }

{\small \noindent[V1568] Giorgio Vasari, \emph{Lives of the artists},
translation by George Bull, Penguin Classics, 1965. }

{\small \noindent[W92a] R.F.C. Walters, Categories and Computer Science,
Carslaw Publications 1991, Cambridge University Press 1992. }

{\small \noindent\lbrack W92b] R.F.C. Walters, An imperative language based on
distributive categories, Mathematical Structures in Computer Science,
2:249--256, 1992. }
\end{document}